\documentclass[12pt]{article}
\usepackage{epsf}
\usepackage{amsmath}
\usepackage{amsthm}

\renewcommand{\epsfsize}[2]{\textwidth}

\theoremstyle{definition}
\newtheorem{ex}{Example}[section]

\title{An Implementation of the Best\-vina-Handel Algorithm for Surface
Homeomorphisms}
\date{}
\author{Peter Brinkmann}

\begin{document}
\maketitle

\begin{abstract}
In \cite{hb2}, Best\-vina and Handel describe an effective algorithm that
determines whether a given homeomorphism of an orientable, possibly punctured
surface is pseudo-Anosov.
We present a software package in Java that realizes
this algorithm for surfaces with one puncture. Moreover, the
package allows the user to define homeomorphisms in terms of Dehn twists, and
in the pseudo-Anosov case it generates images of train tracks in the
sense of Best\-vina-Handel.
\end{abstract}

\section{Introduction}
The fundamental group of a surface $S$ of genus $g$ with one puncture is a free
group $F$ on $2g$ generators. A homeomorphism of $S$ induces
an outer automorphism $\mathcal O$ of $F$, and we can
represent $\mathcal O$ as a homotopy equivalence $f:G\rightarrow G$ of a
finite graph $G\subset S$ homotopy equivalent to $S$.
\par

$f:G\rightarrow G$ is said to be a {\it train track map} if for every $n\geq 1$
and for every edge $e$ of $G$, the restriction of $f^n$ to
the interior of $e$ is an immersion. In
\cite{hb2}, Best\-vina and Handel give an effective algorithm that takes
a homotopy equivalence $f:G\rightarrow G$ representing an
outer automorphism $\mathcal O$ and attempts to find a train track representative
$f':G'\rightarrow G'$ of $\mathcal O$, where $G'$, like $G$, is embedded in and
homotopy equivalent to $S$. If $\mathcal O$ is irreducible\footnote{See
\cite{hb1} for a definition of irreducibility. For our purposes, it is
sufficient to know that an outer automorphism induced by a pseudo-Anosov
homeomorphism of a surface with one puncture will always be irreducible.},
the algorithm will
always succeed. If $\mathcal O$ is reducible, it will either find a train track
representative, or it will conclude that $\mathcal O$ is reducible.
\par

Given a train track representative $f:G\rightarrow G$ of an outer automorphism
induced by a surface homeomorphism $\phi:S\rightarrow S$, Best\-vina and Handel
(see \cite{hb2})
construct a so-called train track\footnote{Thurston's notion of train tracks
is slightly different from the notion of train tracks according to Best\-vina
and Handel. For an exposition of Thurston's theory of surface
homeomorphisms, see \cite{travaux}.}
$\tau$, which can be thought of as being embedded in $S$. 
Using $\tau$, one can effectively decide whether $\phi$ is pseudo-Anosov.
Furthermore, in the pseudo-Anosov case, the following information
can be extracted from $\tau$ and $f$:
\begin{itemize}
\item the growth rate of $\phi$
\item the structure of the stable and unstable foliations of $\phi$, in
particular singular points of the foliations and their indices
\end{itemize}

The software package implements this theory in the case of surfaces of
genus at least two with exactly one puncture. This restriction is motivated
by the fact pseudo-Anosov homeomorphisms of surfaces with one puncture induce
irreducible automorphisms of the fundamental group. This is not true
for surfaces with more than one puncture, and handling this case would
require the implementation of a more complicated algorithm. However,
the theory developed in \cite{hb2} works in full generality (including the
case of closed surfaces, which can be reduced to the case of
punctured surfaces by removing the orbit of a periodic point).

The package consists of three main parts:
\begin{itemize}
\item The first part takes a surface homeomorphism $\phi:S\rightarrow S$
defined by a sequence of
Dehn twists and turns it into a homotopy equivalence of a graph.
\item The second part takes a homotopy equivalence of a graph and either finds
a reduction or a train track representative.
\item The third part constructs a train track $\tau$ from a train track
representative and generates an image of $\tau$ embedded in the
surface $S$.
\end{itemize}
The output of the second and third part combined contain all the
information about $\phi$ listed above. In particular, they decide whether
$\phi$ is pseudo-Anosov.
\par
The package is highly modular, and
the three parts can be used independently. For example, the handling of
Dehn twists has applications beyond the scope of this paper, and the
second part also works for nongeometric outer automorphisms of free groups
(see \cite{hb1}). Moreover, each of the three parts falls into several
functional units, many of which (such as computations and graphics in 
the hyperbolic plane) may be used in other contexts.
\par

The implementation of the train track algorithm is a part of my master's
thesis (\cite{pb1}); the design and implementation of the programs handling
Dehn twists and graphics is a part of my Diplom thesis (\cite{pb2}). It is
my pleasure to express my gratitude to Klaus Johannson and
Werner Ballmann for their help in writing my theses in Knoxville,
Paris and Bonn. I would like to thank Steve Gersten for
encouraging me to write this article. Finally, I would like to thank the
editors and the referees for valuable suggestions.

The software package, written in Java, is available free of charge at
{\tt http://www.math.utah.edu/\~{ }brinkman.} An older version of the
package, written in ANSI-C, is also available. Both versions are
portable and should run on most systems.

\section{Related algorithms and implementations}
There are at three least other implementations of the Best\-vina-Handel
algorithm, each with an emphasis different from the implementation described
here.
\begin{itemize}
\item T.\ White's ''FOLDTOOL'' (\cite{twhite}) is an implementation of the
train track algorithm from \cite{hb1} for free groups. Automorphisms are
entered and displayed as homotopy equivalences of graphs.
\item B.\ Menasco and J.\ Ringland (see \cite{mr})
have implemented the Best\-vina-Handel algorithm in the case of
automorphisms of punctured spheres. Homeomorphisms can be entered as
braid words or as homotopy equivalences of graphs.
Results are displayed as homotopy equivalences of graphs.
\item T.\ Hall's implementation (see \cite{thall}) handles arbitrary punctured
surfaces. Homeomorphisms are entered as homotopy equivalences of graphs,
as braid words, or as horseshoe maps according to Smale.
Results are displayed as homotopy equivalences of graphs.
\end{itemize}

A common characteristic of all implementations is a program realizing some
part of the theory developed in \cite{hb1,hb2}. The main distinguishing
characteristic of the implementation discussed here is that homeomorphisms
of surfaces with one puncture can be entered as compositions of Dehn twists,
and results can be displayed as pictures of graphs embedded in surfaces, which
significantly facilitates the generation of examples as well as the
interpretation of results. Hence, the software described here provides a
powerful yet easy-to-use environment for mathematical experimentation.
\par

Finally, we note that various authors have found independent
approaches to train tracks, e.\ g., M.\ Lustig in \cite{mlhabil}, J.\ Los in
\cite{los} and J.\ Franks/M.\ Misiurewicz in \cite{fm}. In \cite{mlhabil},
train tracks are used to study automorphisms of free groups, while the other
two papers are concerned with homeomorphisms of punctured spheres.

\section{Dehn twists}
\begin{figure}
\centerline{\epsfbox{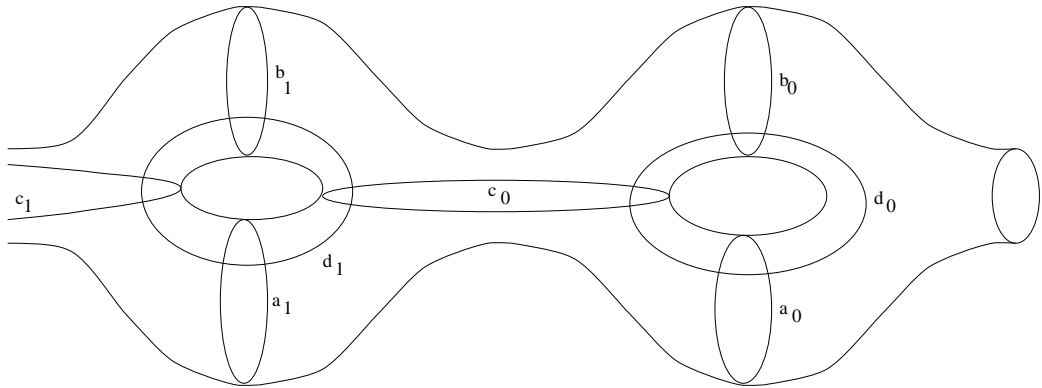}}
\caption{Generators of the mapping class group}\label{mcggen}
\end{figure}

The software package contains a class with two methods for handling Dehn
twists: One of them is extremely easy to use and allows the user to define
surface homeomorphisms
as a composition of Dehn twists with respect to a fixed set of curves
(see figure \ref{mcggen}). The Dehn twists with respect to this set of curves
generate the mapping class group (see \cite{lickorish}). This set of generators
is not minimal; rather, it was chosen with the user's convenience in mind.
\par

The other method for handling Dehn twists removes the restriction to a fixed
set of curves, which results in a slightly more complicated input format. This
method is the part of the package that provides the link between surface
homeomorphisms and homotopy equivalences of graphs; the method described in
the previous paragraph merely generates input for the second one.
\par

When computing Dehn twists, we adopt the following convention: We equip the
surface with an outward pointing normal vector field. When twisting with
respect to a curve $c$, we turn {\it right\footnote{The notion of turning
left or right is defined with respect to the chosen normal vector field.}}
whenever we hit $c$.
\par

\section{Examples}
Figures \ref{nosing}, \ref{onesg}, \ref{minex} and \ref{batman}
were generated by the software package. Each of them
shows a train track belonging to a pseudo-Anosov homeomorphism of a once
punctured surface of genus 2 or 3. The identification pattern on the boundary
of the polygons is given by matching labels of edges intersecting the
boundary, and the puncture corresponds to the vertices of the polygon.
\par
Singularities of the stable or unstable foliation of the pseudo-Anosov
map in question correspond either to the puncture or to shaded areas
containing at least three edges. If a shaded area contains $k\geq 3$ edges,
it gives rise to a singularity of index $1-\frac{k}{2}$. For the proofs of
these statements, see \cite{hb2}.
\par
Since the sum of
the indices of all singularities equals the Euler characteristic of the
surface with the puncture closed, we can compute the index of the singularity
at the puncture, if any.
Moreover, the singularities of the two foliations are fixed points or
periodic points of the pseudo-Anosov homeomorphism in question.
There are more periodic points than just the singularities of the foliations
--- in fact, the set of periodic points of a pseudo-Anosov homeomorphism is
dense, see \cite[expos\'e 9, proposition 18]{travaux}.
\par

In the following examples, $S_g$ is a surface of genus $g$ with one puncture,
and $D_c$ denotes the Dehn twist with respect to a curve $c$, which will
always be one of the curves from figure \ref{mcggen}. All the results in the
following paragraphs were computed by the software package, the only input
being the genus of the surface and a sequence of Dehn twists.
\par

\begin{ex}[Maximal index I]
Consider the map $h:S_2\rightarrow S_2,$
$$ h=D_{a_1}D_{c_0}D_{d_0}D_{a_1}D_{d_1}D_{a_1}.$$
Using the algorithm from \cite{hb2}, the software concludes that $h$ is
a pseudo-Anosov homeomorphism with growth rate $\lambda\approx 1.722084$.
A train track for $h$ is shown in figure
\ref{nosing}. None of the shaded areas gives rise to a singularity of the
stable or unstable foliation, so the puncture is the only singularity, and
its index is $-2$.

\begin{figure}
\centerline{\epsfbox{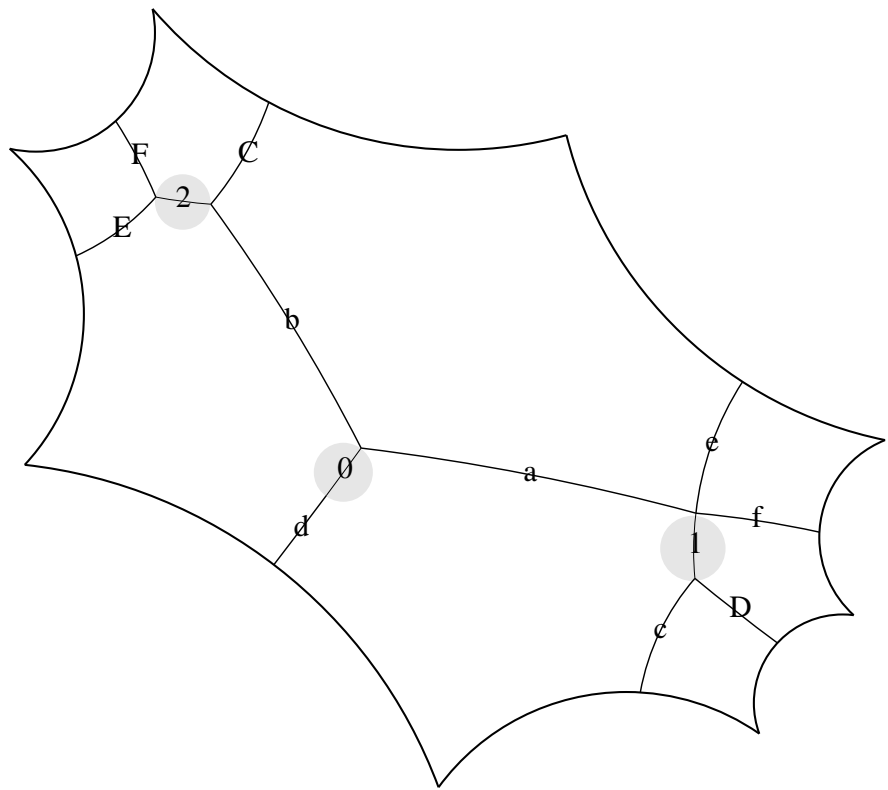}}
\caption{Maximal index I}\label{nosing}
\end{figure}
\end{ex}

\begin{ex}[Maximal index II]
Let $h:S_2\rightarrow S_2$ be given by
$$h=D^{-1}_{a_1}D_{d_1}D^{-1}_{c_0}D_{d_0}.$$
$h$ is a pseudo-Anosov
homeomorphism with growth rate $\lambda\approx 4.390257$. Figure \ref{onesg}
shows the corresponding train track.
The unique shaded area in figure \ref{onesg} contains six edges, so it gives
rise to a singularity $p$ of index $-2$. We conclude that there is no
singularity at the puncture.

\begin{figure}
\centerline{\epsfbox{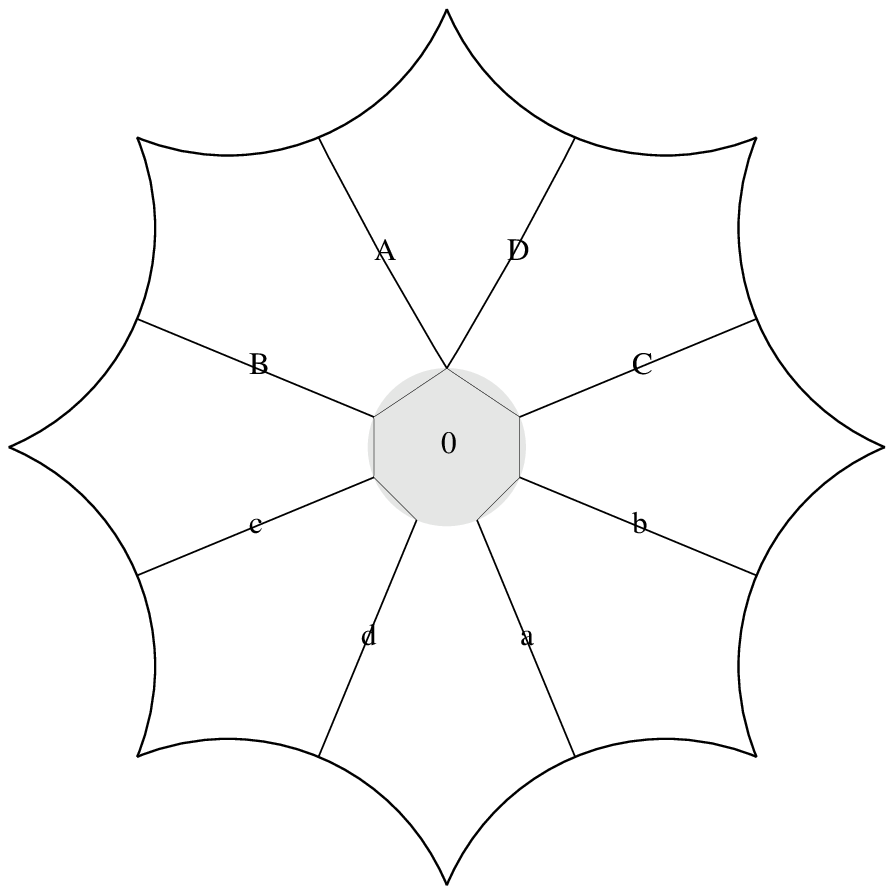}}
\caption{Maximal index II}\label{onesg}
\end{figure}
\end{ex}

\begin{ex}[Minimal index]
Let the homeomorphism $h:S_2\rightarrow S_2$ be given by
$$h=D_{a_0}D^{-1}_{c_0}D_{d_0}D^{-1}_{d_1}.$$
$h$ is a pseudo-Anosov
homeomorphism with growth rate $\lambda\approx 2.015357$. Figure \ref{minex}
shows the corresponding train track.
The shaded areas labeled $0,1,3,4$ give rise to singularities of index
$-\frac{1}{2}$, which shows that there is no singularity at the puncture.
The singularities $0$ and $4$ as well as $1$ and $3$ are exchanged by $h$.

\begin{figure}
\centerline{\epsfbox{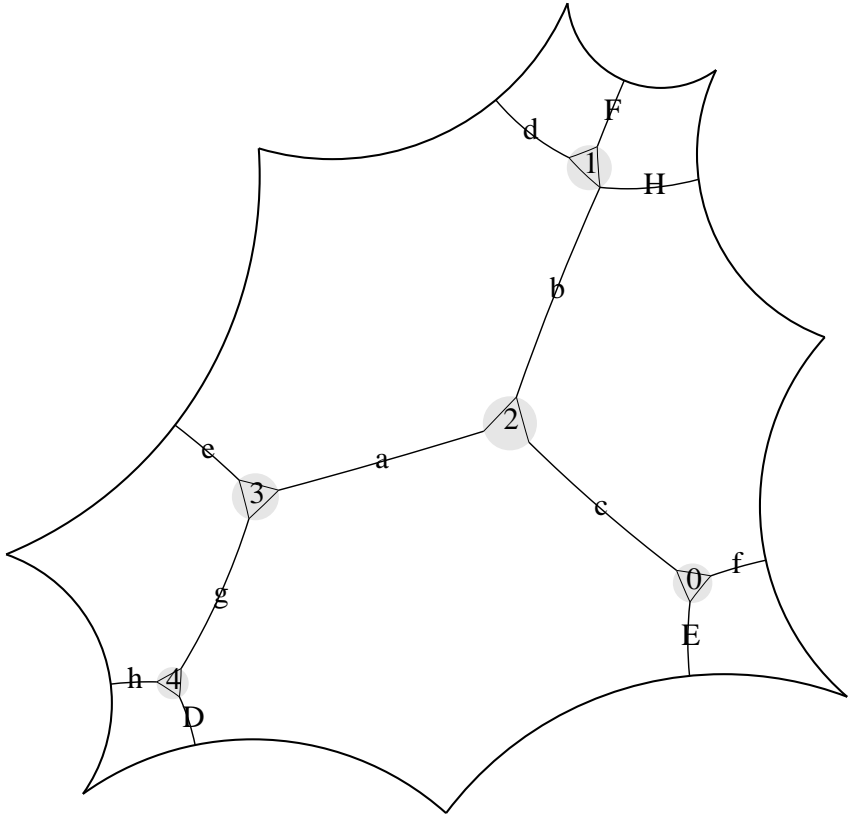}}
\caption{Minimal index}\label{minex}
\end{figure}
\end{ex}

\begin{ex}[Genus 3]
Let $h:S_3\rightarrow S_3$ be given by
$$h=D_{d_0}D_{c_0}D_{d_1}D_{c_1}D_{d_2}D^{-1}_{c_2}.$$
$h$ is a pseudo-Anosov
homeomorphism with growth rate $\lambda\approx 2.042491$. Figure \ref{batman}
shows the corresponding train track.
The shaded areas labeled $0, 2$ give rise to singularities of index $-2$, and
they are exchanged by $h$. There is no singularity at the puncture.

\begin{figure}
\centerline{\epsfbox{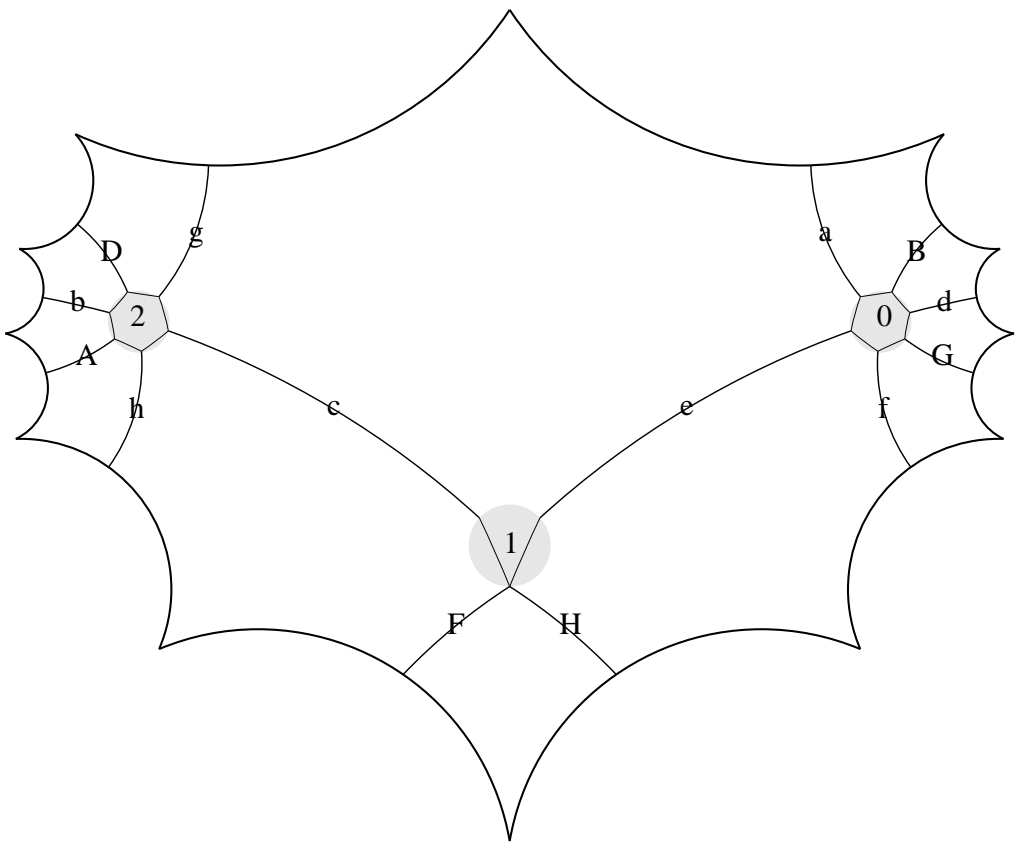}}
\caption{Example of genus 3}\label{batman}
\end{figure}
\end{ex}

\begin{ex}[A reducible example]
Finally, consider $h:S_2\rightarrow S_2$ defined by
$$h=D_{d_0}D_{c_0}D_{d_1}.$$
$h$ is reducible since the complement of
the curves $d_0$, $c_0$, and $d_1$ is not a (punctured) disc, and in fact
the software reaches the same conclusion.
\end{ex}

\section{Implementation}
The complete online documentation of the software package, including a user
manual and the source code, can be found at
{\tt http://www.math.utah.edu/ \~{ }brinkman,} so at this
point we restrict ourselves to a brief discussion of the main implementation
issues. For the most part, we take the point of view of
mathematics rather than that of computer science.
\par

\subsection{Encoding of embeddings}
For the rest of this discussion, it will be advantageous to think of punctures
as being distinguished points of closed surfaces.
Given a closed surface $S$ with a distinguished point $p$ and a finite graph
$G\subset S$ homotopy equivalent to $S-\{p\}$,
we need an efficient way of encoding the embedding of $G$
in $S$. To this end, consider a loop $\rho'$ around $p$. $\rho'$ is homotopic
to a closed edge path $\rho$ in $G$ that crosses every edge of $G$ twice,
once for each direction (assuming that $G$ has no vertices of valence one). 
Conversely, given $G$ and $\rho$, we can reconstruct $S$: We simply take a
polygon $P$ with $2n$ sides, where $n$ is the number of edges of $G$, and
interpret $\rho$ as an identification pattern on the boundary of $P$.
Moreover, we can triangulate $P$ (and hence $S$) by fixing a point $p$
in the interior of $P$ and connecting $p$ to all the vertices in the
boundary of $P$. Hence, we see that
$G$ and $\rho$ give us an efficient way of encoding the embedding of $G$ in
$S$ along with a triangulation of $S$.
\par

\subsection{Finding a metric}
Now, given a triangulation $\tau$ of $S$, we want to find a hyperbolic metric
on $S$ with the property that the edges of $\tau$ are geodesic segments. There
are various ways to accomplish this (see \cite{colin}); our method of choice
is a special case of Thurston's {\it circle packing} (see \cite{notes}):
Given a surface
$S$ with a hyperbolic metric $\mu$ and with a triangulation $\tau$ whose
edges are geodesic segments, there is a collection of circles centered at the
vertices of $\tau$ such that no two circles intersect transversally and two
vertices of $\tau$ are connected by an edge if and only if their corresponding
circles are tangent. For each triangulation, there exists exactly one such
set of circles, and their radii
can be computed numerically. Moreover, they uniquely determine $\mu$.
Hence, circle packing gives us an effective way
of drawing $S$ as a polygon (with identifications on the boundary) in
the hyperbolic plane.

\subsection{Philosophy}
The package takes advantage of many features of the object-oriented paradigm,
such as data encapsulation and reusability. For example, the class that
implements maps of graphs does not allow direct access to its contents; the
other parts of the package operate on such maps through a small and well
defined set of methods, which results in ease of maintenance and great
flexibility.

At this time, one general drawback of Java is that many web browsers
use faulty implementations of the Java libraries, which may cause problems
when running the software described in this article. However, such problems
do not occur when the software is run by the Sun appletviewer, which is
available for most systems.

The mathematical part of the package consists of 16 classes,
reflecting increasing levels
of specialization. Some of them, like the implementation of the train track
algorithm from \cite{hb1}, will only be used in the context of this package.
Others, like the collection of methods for computations and drawings in the
hyperbolic plane, have been designed with other uses in mind. In fact, the
package presented here does not even use all the methods defined in this
collection.
\par

Finally, the classes and methods handling maps of graphs may be
useful beyond the context of this article. For example, the author has
already used them for a tentative implementation of some of the algorithms
in \cite{stallings}.
\par

\bibliographystyle{is-alpha}
\bibliography{all}

\begin{thebibliography}{CdV91}
\ifx \showCODEN  \undefined \def \showCODEN #1{CODEN #1}  \fi
\ifx \showISBN   \undefined \def \showISBN  #1{ISBN #1}   \fi
\ifx \showISSN   \undefined \def \showISSN  #1{ISSN #1}   \fi
\ifx \showLCCN   \undefined \def \showLCCN  #1{LCCN #1}   \fi
\ifx \showPRICE  \undefined \def \showPRICE #1{#1}        \fi
\ifx \showURL    \undefined \def \showURL {URL }          \fi
\ifx \path       \undefined \input path.sty               \fi
\ifx \ifshowURL \undefined
     \newif \ifshowURL
     \showURLtrue
\fi

\bibitem[BH92]{hb1}
Mladen Bestvina and Michael Handel.
\newblock Train tracks and automorphisms of free groups.
\newblock {\em Annals of Mathematics}, 135:\penalty0 1--51, 1992.

\bibitem[BH95]{hb2}
Mladen Bestvina and Michael Handel.
\newblock Train tracks for surface homeomorphisms.
\newblock {\em Topology}, 34\penalty0 (1):\penalty0 109--140, 1995.

\bibitem[Bri95]{pb1}
Peter Brinkmann.
\newblock {Pseudo-Anosov Automorphisms of Free Groups}.
\newblock Master's thesis, University of Tennessee, August 1995.
\newblock \ifshowURL {\showURL
  \path|http://www.math.utah.edu/~brinkman/train|}\fi.

\bibitem[Bri96]{pb2}
Peter Brinkmann.
\newblock {Ein algorithmischer Zugang zur Klassifikation von
  Fl\"achenhom\"oomorphismen}.
\newblock Diplomarbeit, University of Bonn, December 1996.
\newblock \ifshowURL {\showURL
  \path|http://www.math.utah.edu/~brinkman/train|}\fi.

\bibitem[CdV91]{colin}
Yves Colin~de Verdi\`ere.
\newblock Comment rendre g\'eod\'esique une triangulation d'une surface?
\newblock {\em L'Enseignement Math\'ematique}, 37:\penalty0 201--212, 1991.

\bibitem[FLP79]{travaux}
A.~Fathi, F.~Laudenbach, and V.~Poenaru.
\newblock {\em Travaux de Thurston sur les Surfaces}, volume 66 - 67 of {\em
  Ast\'erisque}.
\newblock 1979.

\bibitem[FM93]{fm}
J.~Franks and M.~Misiurewicz.
\newblock {Cycles for disk homeomorphisms and thick trees}.
\newblock {\em Cont.\ Math.}, 152:\penalty0 69--139, 1993.

\bibitem[Hal96]{thall}
T.~Hall.
\newblock Train tracks of surface homeomorphisms: Computer software.
\newblock 1996.
\newblock \ifshowURL {\showURL
  \path|http://www.liv.ac.uk/PureMaths/members/T_Hall.html|}\fi.

\bibitem[Lic64]{lickorish}
W.~Lickorish.
\newblock {A finite set of generators for the homeotopy group of a
  $2$-manifold}.
\newblock {\em Proc. Cambridge Philos. Soc.}, 60:\penalty0 769--778, 1964.

\bibitem[Los93]{los}
J.~Los.
\newblock {Pseudo-Anosov maps and invariant train tracks in the disc: a finite
  algorithm}.
\newblock {\em Proc. London Math. Soc.}, 66\penalty0 (3):\penalty0 400--430,
  1993.

\bibitem[Lus92]{mlhabil}
Martin Lustig.
\newblock {Automorphismen von freien Gruppen}.
\newblock Ha\-bi\-li\-ta\-tions\-schrift, 1992.

\bibitem[MR96]{mr}
W.~Menasco and J.~Ringland.
\newblock {BH2.1: An Implementation of the Bestvina-Handel Algorithm}.
\newblock 1996.
\newblock \ifshowURL {\showURL
  \path|http://emerald.math.buffalo.edu/~ringland/jade/BH/|}\fi.

\bibitem[Sta83]{stallings}
J.R. Stallings.
\newblock {Topology of finite graphs}.
\newblock {\em Inv. Math.}, 71:\penalty0 551--565, 1983.

\bibitem[Thu78]{notes}
William~P. Thurston.
\newblock The geometry and topology of three-manifolds.
\newblock Princeton notes, 1978.

\bibitem[Whi90]{twhite}
T.~White.
\newblock {FOLDTOOL}.
\newblock software package, 1990.

\end{thebibliography}
\par

{\noindent \sc Department of Mathematics, University of Utah\\
Salt Lake City, UT 84112, USA\\}
\par
{\noindent \it E-mail:} brinkman\@@math.utah.edu

\end{document}